\documentclass[english,11pt]{article}
\usepackage[latin1]{inputenc}
\usepackage{amsmath,amsthm,amssymb,babel}

\newtheorem{teo}{Theorem}

\newtheorem{lemma}{Lemma}

\def\proof{{\it Proof.}\ }
\def\endproof{\hfill $\Box$\par\vskip3mm}

\def\eq#1{(\ref{#1})}

\def\neweq#1{\begin{equation}\label{#1}}
\def\endeq{\end{equation}}

\def\phi{\varphi}
\def\RR{{\mathbb R} }

\def\di{\displaystyle}

\date{}

\title{\sc Continuous spectrum for a class of nonhomogeneous differential operators \thanks{
Correspondence address: Vicen\c{t}iu R\u{a}dulescu, Department of
Mathematics, University of Craiova,  200585 Craiova, Romania. E-mail:
{\tt
vicentiu.radulescu@math.cnrs.fr}}}
\author{\sc Mihai Mih\u ailescu$\,^a$ and Vicen\c{t}iu R\u{a}dulescu$\,^{a,b}$\\
\small
$^a\,$Department of Mathematics, University of Craiova,  200585 Craiova,
Romania\\
\small $^b\,$Institute of Mathematics ``Simion Stoilow" of the Romanian Academy,\\
\small P.O. Box 1-764, 014700 Bucharest, Romania\\
\small
E-mail addresses: {\tt mmihailes@yahoo.com}\qquad {\tt
vicentiu.radulescu@math.cnrs.fr}}

\begin{document}
\baselineskip16pt \maketitle \noindent{\small{\sc Abstract}. We
study the boundary value problem $-{\rm div}((|\nabla
u|^{p_1(x)-2}+|\nabla u|^{p_2(x)-2})\nabla
u)=\lambda|u|^{q(x)-2}u$ in $\Omega$, $u=0$ on $\partial\Omega$,
where $\Omega$ is a bounded domain in $\RR^N$ with smooth
boundary, $\lambda$ is a positive real number, and the continuous
functions $p_1$, $p_2$, and $q$ satisfy $1<p_2(x)<q(x)<p_1(x)<N$
and $\max_{y\in\overline\Omega}q(y)<\frac{N p_2(x)}{N-p_2(x)}$ for
any $x\in\overline\Omega$. The main result of this paper
establishes the existence of two positive constants $\lambda_0$
and $\lambda_1$ with $\lambda_0\leq\lambda_1$ such that any
$\lambda\in[\lambda_1,\infty)$ is an eigenvalue,  while
any $\lambda\in(0,\lambda_0)$ is not an eigenvalue of the above problem.
}

\section{Introduction and preliminary results}
In this paper we are concerned with the study of the eigenvalue problem
\begin{equation}\label{1}
\left\{\begin{array}{lll}
-{\rm div}((|\nabla u|^{p_1(x)-2}+|\nabla u|^{p_2(x)-2})\nabla u)=\lambda|u|^{q(x)-2}u, &\mbox{for}&
x\in\Omega\\
u=0, &\mbox{for}& x\in\partial\Omega\,,
\end{array}\right.
\end{equation}
where $\Omega\subset\RR^N$ ($N\geq 3$) is a bounded domain with
smooth boundary, $\lambda>0$ is a real number, and $p_1$, $p_2$, $q$ are
continuous functions on $\overline\Omega$.

The study of eigenvalue problems involving operators with variable
exponents growth conditions has captured a special attention in
the last few years. This is in keeping with the fact that
operators which arise in such kind of problems, like the
$p(x)$-Laplace operator (i.e., ${\rm div}(|\nabla
u|^{p(x)-2}\nabla u)$, where $p(x)$ is a continuous positive
function),  are not homogeneous and thus,  a large number of
techniques which can be applied in the homogeneous case (when
$p(x)$ is a positive constant) fail  in this new setting. A
typical example is the Lagrange multiplier theorem, which does not
apply to the eigenvalue problem
\begin{equation}\label{P}
\left\{\begin{array}{lll}
-{\rm div}(|\nabla u|^{p(x)-2}\nabla u)=\lambda|u|^{q(x)-2}u, &\mbox{for}&
x\in\Omega\\
u=0, &\mbox{for}& x\in\partial\Omega\,,
\end{array}\right.
\end{equation}
 where
$\Omega\subset\RR^N$  is a bounded domain. This is due to the fact
that the associated Rayleigh quotient is not homogeneous, provided
both $p$ and $q$ are not constant.

On the other hand, problems like \eq{P} have been largely
considered in the literature in the recent years. We  give in what
follows a concise but complete image of the actual stage of
research on this topic.
\medskip

$\bullet$ In the case when $p(x)=q(x)$ on $\overline\Omega$,  Fan, Zhang and Zhao \cite{FZZ} established  the existence of
infinitely many eigenvalues for problem \eq{P} by using an argument based on the Ljusternik-Schnirelmann
critical point theory. Denoting by $\Lambda$ the set of all
nonnegative eigenvalues, Fan, Zhang and Zhao showed that $\Lambda$ is discrete, $\sup\Lambda=+\infty$ and they pointed out that
 only under special conditions, which are somehow connected with a kind of monotony of the function $p(x)$, we have
$\inf\Lambda>0$ (this is in contrast with the case when $p(x)$ is a
constant; then, we always have  $\inf\Lambda>0$).
\medskip

$\bullet$ In the case when $\min_{x\in\overline\Omega}q(x)<\min_{x\in\overline\Omega}p(x)$ and
$q(x)$ has a subcritical growth Mih\u ailescu and R\u adulescu \cite{mihradproc} used the Ekeland's variational
principle in order to prove the existence of a continuous family of eigenvalues  which lies in a neighborhood of
the origin.
\medskip

$\bullet$ In the case when  $\max_{x\in\overline\Omega}p(x)<\min_{x\in\overline\Omega}q(x)$ and
$q(x)$ has a subcritical growth a mountain-pass argument, similar with those used by Fan and Zhang in
the proof of Theorem 4.7 in \cite{FZh}, can be applied in order to show that any $\lambda>0$ is an
eigenvalue of problem \eq{P}.
\medskip

$\bullet$ In the case when
$\max_{x\in\overline\Omega}q(x)<\min_{x\in\overline\Omega}p(x)$ it
can be proved that the energy functional associated to problem
\eq{P} has a nontrivial minimum for any positive $\lambda$ large
enough (see Theorem 4.7 in \cite{FZh}). Clearly, in this case the
result in \cite{mihradproc} can be also applied. Consequently, in
this situation there exist two positive constants $\lambda^\star$
and $\lambda^{\star\star}$ such that any
$\lambda\in(0,\lambda^\star)\cup (\lambda^{\star\star},\infty)$ is
an eigenvalue of problem \eq{P}.

\smallskip In this paper we study problem \eq{1} under the
following assumptions:
\begin{equation}\label{2}
1<p_2(x)<\min_{y\in\overline\Omega}q(y)\leq\max_{y\in\overline\Omega}q(y)<p_1(x)<N,\;\;\;\forall\;x\in\overline\Omega
\end{equation}
and
\begin{equation}\label{3}
\max_{y\in\overline\Omega}q(y)<\frac{N
p_2(x)}{N-p_2(x)},\;\;\;\forall\;x\in\overline\Omega\, .
\end{equation}
Thus, the case considered here is different from all the cases
studied before. In this new situation we will show the existence
of two positive constants $\lambda_0$ and $\lambda_1$ with
$\lambda_0\leq\lambda_1$ such that any
$\lambda\in[\lambda_1,\infty)$ is an eigenvalue of problem \eq{1}
while any $\lambda\in(0,\lambda_0)$ is not an eigenvalue of
problem \eq{1}. An important consequence of our study is that,
under hypotheses \eq{2} and \eq{3}, we have
$$\inf\limits_{u\in W_0^{1,p_1(x)}(\Omega)\setminus\{0\}}
\di\frac{\di\int_\Omega\frac{1}{p_1(x)}|\nabla
u|^{p_1(x)}\;dx+\di\int_\Omega\frac{1}{p_2(x)}|\nabla
u|^{p_2(x)}\;dx}
{\di\int_\Omega\frac{1}{q(x)}|u|^{q(x)}\;dx}>0\,.$$ That fact is
proved by using the Lagrange Multiplier Theorem. The absence of
homogeneity will be balanced  by the fact that assumptions \eq{2}
and \eq{3} yield
$$\lim_{\|u\|_{p_1(x)}\rightarrow 0}\di\frac{\di\int_\Omega\frac{1}{p_1(x)}|\nabla
u|^{p_1(x)}\;dx+\di\int_\Omega\frac{1}{p_2(x)}|\nabla
u|^{p_2(x)}\;dx}
{\di\int_\Omega\frac{1}{q(x)}|u|^{q(x)}\;dx}=\infty$$ and
$$\lim_{\|u\|_{p_1(x)}\rightarrow\infty}\di\frac{\di\int_\Omega\frac{1}{p_1(x)}|\nabla
u|^{p_1(x)}\;dx+\di\int_\Omega\frac{1}{p_2(x)}|\nabla
u|^{p_2(x)}\;dx}
{\di\int_\Omega\frac{1}{q(x)}|u|^{q(x)}\;dx}=\infty\,,$$ where
$\|\cdot\|_{p_1(x)}$ stands for the norm in the variable exponent
Sobolev space $W_0^{1,p_1(x)}(\Omega)$.

\smallskip
We start with some preliminary basic results on the theory of
 Lebesgue--Sobolev spaces with variable exponent. For more details we refer to the book by Musielak
\cite{M} and the papers by Edmunds et al. \cite{edm, edm2, edm3},
Kovacik and  R\'akosn\'{\i}k \cite{KR}, Mih\u ailescu and R\u
adulescu \cite{RoyalSoc,mihradjmaa}, and Samko and Vakulov
\cite{samko}.

Set
$$C_+(\overline\Omega)=\{h;\;h\in C(\overline\Omega),\;h(x)>1\;{\rm
for}\;
{\rm all}\;x\in\overline\Omega\}.$$
For any $h\in C_+(\overline\Omega)$ we define
$$h^+=\sup_{x\in\Omega}h(x)\qquad\mbox{and}\qquad h^-=
\inf_{x\in\Omega}h(x).$$ For any $p\in C_+(\overline\Omega)$, we
define the variable exponent Lebesgue space
$$L^{p(x)}(\Omega)=\left\{u;\ u\ \mbox{is a
 measurable real-valued function such that }
\int_\Omega|u(x)|^{p(x)}\;dx<\infty\right\}.$$ We define on this
space the  {\it Luxemburg norm} by
$$|u|_{p(x)}=\inf\left\{\mu>0;\;\int_\Omega\left|
\frac{u(x)}{\mu}\right|^{p(x)}\;dx\leq 1\right\}.$$

Let $L^{p^{'}(x)}(\Omega)$ denote the conjugate space of
$L^{p(x)}(\Omega)$, where $1/p(x)+1/p^{'}(x)=1$. For any $u\in
L^{p(x)}(\Omega)$ and $v\in L^{p^{'}(x)}(\Omega)$ the H\"older
type inequality
\begin{equation}\label{Hol}
\left|\int_\Omega uv\;dx\right|\leq\left(\frac{1}{p^-}+
\frac{1}{{p^{'}}^-}\right)|u|_{p(x)}|v|_{p^{'}(x)}
\end{equation}
holds true.

An important role in manipulating the generalized Lebesgue-Sobolev
spaces
is played by the {\it modular} of the $L^{p(x)}(\Omega)$ space, which
is
the mapping
 $\rho_{p(x)}:L^{p(x)}(\Omega)\rightarrow\RR$ defined by
$$\rho_{p(x)}(u)=\int_\Omega|u|^{p(x)}\;dx.$$
If $(u_n)$, $u\in L^{p(x)}(\Omega)$ then the following relations
hold true
\begin{equation}\label{L4}
|u|_{p(x)}>1\;\;\;\Rightarrow\;\;\;|u|_{p(x)}^{p^-}\leq\rho_{p(x)}(u)
\leq|u|_{p(x)}^{p^+}
\end{equation}
\begin{equation}\label{L5}
|u|_{p(x)}<1\;\;\;\Rightarrow\;\;\;|u|_{p(x)}^{p^+}\leq
\rho_{p(x)}(u)\leq|u|_{p(x)}^{p^-}
\end{equation}
\begin{equation}\label{L6}
|u_n-u|_{p(x)}\rightarrow 0\;\;\;\Leftrightarrow\;\;\;\rho_{p(x)}
(u_n-u)\rightarrow 0.
\end{equation}

Next, we define $W_0^{1,p(x)}(\Omega)$ as the closure of
$C_0^\infty(\Omega)$ under the norm
$$\|u\|_{p_(x)}=|\nabla u|_{p(x)}.$$
The space $W_0^{1,p(x)}(\Omega)$ is a separable and reflexive
Banach space. We note that if $s\in C_+(\overline\Omega)$ and
$s(x)<p^\star(x)$ for all $x\in\overline\Omega$ then the embedding
$W_0^{1,p(x)}(\Omega)\hookrightarrow L^{s(x)}(\Omega)$ is compact
and continuous, where $p^\star(x)=\frac{Np(x)}{N-p(x)}$ if
$p(x)<N$ or $p^\star(x)=+\infty$ if $p(x)\geq N$.

For applications of Sobolev spaces with variable exponent we refer
to Acerbi and Mingione \cite{acerbi2}, Chen, Levine and Rao
\cite{chen}, Diening \cite{D},  Halsey \cite{hal}, Ruzicka
\cite{R}, and Zhikov \cite{Z1}).

\section{The main result}
Since $p_2(x)<p_1(x)$ for any $x\in\overline\Omega$ it follows that
$W_0^{1,p_1(x)}(\Omega)$ is continuously embedded in
$W_0^{1,p_2(x)}(\Omega)$. Thus, a solution for a problem of type
\eq{1} will be sought in the variable exponent space
$W_0^{1,p_1(x)}(\Omega)$.

We say that $\lambda\in\RR$ is an {\it eigenvalue} of problem \eq{1} if
there exists $u\in W_0^{1,p_1(x)}(\Omega)
\setminus\{0\}$ such that
$$\int_{\Omega}(|\nabla u|^{p_1(x)-2}+|\nabla u|^{p_2(x)-2})\nabla u\nabla v\;dx-\lambda\int_{\Omega}
|u|^{q(x)-2}uv\;dx=0\,,$$ for all $v\in W_0^{1,p_1(x)}(\Omega)$.
We point out that if $\lambda$ is an eigenvalue of problem \eq{1}
then the corresponding eigenfunction $u\in
W_0^{1,p_1(x)}(\Omega)\setminus\{0\}$ is a {\it weak solution} of
problem \eq{1}.

Define
$$\lambda_1:=\inf\limits_{u\in W_0^{1,p_1(x)}(\Omega)\setminus\{0\}}
\di\frac{\di\int_\Omega\frac{1}{p_1(x)}|\nabla
u|^{p_1(x)}\;dx+\di\int_\Omega\frac{1}{p_2(x)}|\nabla
u|^{p_2(x)}\;dx} {\di\int_\Omega\frac{1}{q(x)}|u|^{q(x)}\;dx}\,.$$

Our main result is given by the following theorem.

\begin{teo}\label{t1}
Assume that conditions \eq{2} and \eq{3} are fulfilled. Then $\lambda_1>0$. Moreover, any
$\lambda\in[\lambda_1,\infty)$ is an eigenvalue of problem \eq{1}. Furthermore, there
exists a positive constant $\lambda_0$ such that $\lambda_0\leq\lambda_1$ and any $\lambda\in(0,\lambda_0)$ is not
an eigenvalue of problem \eq{1}.
\end{teo}

\proof Let $E$ denote the generalized Sobolev space
$W_0^{1,p_1(x)}(\Omega)$. We denote by $\|\cdot\|$ the norm on
$W_0^{1,p_1(x)}(\Omega)$ and by $\|\cdot\|_1$ the norm on
$W_0^{1,p_2(x)}(\Omega)$.

Define the functionals $J$, $I$, $J_1$, $I_1:E\rightarrow\RR$ by
$$J(u)=\int_{\Omega}\frac{1}{p_1(x)}|\nabla u|^{p_1(x)}\;dx+\int_{\Omega}\frac{1}{p_2(x)}|\nabla u|^{p_2(x)}\;dx,$$
$$I(u)=\int_\Omega\frac{1}{q(x)}|u|^{q(x)}\;dx,$$
$$J_1(u)=\int_{\Omega}|\nabla u|^{p_1(x)}\;dx+\int_{\Omega}|\nabla u|^{p_2(x)}\;dx,$$
$$I_1(u)=\int_\Omega|u|^{q(x)}\;dx.$$
Standard arguments imply that $J,\;I\in C^1(E,\RR)$ and for all
$u,\;v\in E$,
$$\langle J^{'}(u),v\rangle=\int_{\Omega}(|\nabla u|^{p_1(x)-2}+|\nabla u|^{p_2(x)-2})
\nabla u\nabla v\;dx,$$
$$\langle I^{'}(u),v\rangle=\int_\Omega|u|^{q(x)-2}uv\;dx.$$

We split the proof of Theorem \ref{t1} into four steps.

\noindent$\bullet$ {\sc Step 1.} We show that $\lambda_1>0$.

Since for any $x\in\overline\Omega$ we have $p_1(x)>q^+\geq
q(x)\geq q^->p_2(x)$ we deduce that for any $u\in E$,
$$2(|\nabla u(x)|^{p_1(x)}+|\nabla u(x)|^{p_2(x)})\geq|\nabla u(x)|^{q^+}+|\nabla u(x)|^{q^-}$$
and
$$|u(x)|^{q^+}+|u(x)|^{q^-}\geq|u(x)|^{q(x)}.$$
Integrating the above inequalities we find
\begin{equation}\label{4}
2\int_\Omega(|\nabla u|^{p_1(x)}+|\nabla
u|^{p_2(x)})\;dx\geq\int_\Omega(|\nabla u|^{q^+}+|\nabla
u|^{q^-})\;dx, \qquad\forall\;u\in E
\end{equation}
and
\begin{equation}\label{4prim}
\int_\Omega(|u|^{q^+}+|u|^{q^-})\;dx\geq\int_\Omega|u|^{q(x)}\;dx,\qquad\forall\;u\in
E.
\end{equation}
By Sobolev embeddings, there exist positive constants
$\lambda_{q^+}$ and $\lambda_{q^-}$ such that
\begin{equation}\label{5}
\int_\Omega|\nabla
u|^{q^+}\;dx\geq\lambda_{q^+}\int_\Omega|u|^{q^+}\;dx,\qquad\forall\;u\in
W_0^{1,q^+}(\Omega)
\end{equation}
and
\begin{equation}\label{6}
\int_\Omega|\nabla
u|^{q^-}\;dx\geq\lambda_{q^-}\int_\Omega|u|^{q^-}\;dx,\qquad\forall\;u\in
W_0^{1,q^-}(\Omega).
\end{equation}
Using again the fact that $q^-\leq q^+<p_1(x)$ for any
$x\in\overline\Omega$ we deduce that $E$ is continuously embedded in
$W_0^{1,q^+}(\Omega)$ and in $W_0^{1,q^-}(\Omega)$. Thus,
inequalities \eq{5} and \eq{6} hold true for any $u\in E$.

Using inequalities \eq{5}, \eq{6} and \eq{4prim} it is clear that there exists a positive constant $\mu$ such that
\begin{equation}\label{7}
\int_\Omega(|\nabla u|^{q^+}+|\nabla
u|^{q^-})\;dx\geq\mu\int_\Omega|u|^{q(x)}\;dx,\qquad\forall\;u\in
E.
\end{equation}
Next, inequalities \eq{7} and \eq{4} yield
\begin{equation}\label{8}
\int_\Omega(|\nabla u|^{p_1(x)}+|\nabla
u|^{p_2(x)})\;dx\geq\frac{\mu}{2}\int_\Omega|u|^{q(x)}\;dx,\qquad\forall\;u\in
E.
\end{equation}
By relation \eq{8} we deduce that
\begin{equation}\label{lambda0}
\lambda_0=\inf_{v\in E\setminus\{0\}}\frac{J_1(v)}{I_1(v)}>0
\end{equation}
and thus,
\begin{equation}\label{00}
J_1(u)\geq\lambda_0 I_1(u),\;\;\;\forall\;u\in E.
\end{equation}
The above inequality yields
\begin{equation}\label{9}
p_1^+\cdot J(u)\geq J_1(u)\geq\lambda_0 I_1(u)\geq\lambda_0 I(u)\;\;\;\forall\;u\in E.
\end{equation}
The last inequality assures that $\lambda_1>0$ and thus, step 1 is
verified.
\medskip

\noindent$\bullet$ {\sc Step 2.} We show that $\lambda_1$ is an
eigenvalue of problem \eq{1}.

\begin{lemma}\label{l1}
The following relations hold true:
\begin{equation}\label{10}
\lim_{\|u\|\rightarrow\infty}\frac{J(u)}{I(u)}=\infty
\end{equation}
and
\begin{equation}\label{11}
\lim_{\|u\|\rightarrow 0}\frac{J(u)}{I(u)}=\infty.
\end{equation}
\end{lemma}
\proof Since $E$ is continuously embedded in $L^{q^\pm}(\Omega)$ it
follows that there exist two positive constants $c_1$ and $c_2$ such
that
\begin{equation}\label{12}
\|u\|\geq c_1\cdot|u|_{q^+},\;\;\;\forall\;u\in E
\end{equation}
and
\begin{equation}\label{13}
\|u\|\geq c_2\cdot|u|_{q^-},\;\;\;\forall\;u\in E.
\end{equation}
For any $u\in E$ with $\|u\|>1$ by relations \eq{L4}, \eq{4prim}, \eq{12}, \eq{13} we infer
$$\frac{J(u)}{I(u)}\geq\frac{\di\frac{\|u\|^{p_1^-}}{p_1^+}}{\di\frac{|u|_{q^+}^{q^+}+|u|_{q^-}^{q^-}}{q^-}}
\geq\frac{\di\frac{\|u\|^{p_1^-}}{p_1^+}}{\di\frac{c_1^{-q^+}\|u\|^{q^+}+c_2^{-q^-}\|u\|^{q^-}}{q^-}}.$$
Since $p_1^->q^+\geq q^-$, passing to the limit as $\|u\|\rightarrow\infty$ in the above inequality we deduce that
relation \eq{10} holds true.

Next, let us remark that since $p_1(x)>p_2(x)$ for any
$x\in\overline\Omega$, the space $W_0^{1,p_1(x)}(\Omega)$ is
continuously embedded in $W_0^{1,p_2(x)}(\Omega)$. Thus, if
$\|u\|\rightarrow 0$ then $\|u\|_1\rightarrow 0$.

The above remarks enable us to affirm that for any $u\in E$ with $\|u\|<1$ small enough we have $\|u\|_1<1$.

On the other hand, since \eq{3} holds true we deduce that $W_0^{1,p_2(x)}(\Omega)$ is continuously
embedded in $L^{q^\pm}(\Omega)$.
 It follows that there exist two positive constants
$d_1$ and $d_2$ such that
\begin{equation}\label{12prim}
\|u\|_1\geq d_1\cdot|u|_{q^+},\;\;\;\forall\;u\in W_0^{1,p_2(x)}(\Omega)
\end{equation}
and
\begin{equation}\label{13prim}
\|u\|_1\geq d_2\cdot|u|_{q^-},\;\;\;\forall\;u\in W_0^{1,p_2(x)}(\Omega).
\end{equation}
Thus, for any $u\in E$ with $\|u\|<1$ small enough, relations \eq{L5}, \eq{4prim}, \eq{12prim}, \eq{13prim} imply
$$\frac{J(u)}{I(u)}\geq\frac{\di\frac{\int_\Omega|\nabla u|^{p_2(x)}\;dx}{p_2^+}}{\di\frac{|u|_{q^+}^{q^+}+|u|_{q^-}^{q^-}}{q^-}}
\geq\frac{\di\frac{\|u\|_1^{p_2^+}}{p_2^+}}{\di\frac{d_1^{-q^+}\|u\|_1^{q^+}+d_2^{-q^-}\|u\|_1^{q^-}}{q^-}}.$$
Since $p_2^+<q^-\leq q^+$, passing to the limit as
$\|u\|\rightarrow 0$ (and thus, $\|u\|_1\rightarrow 0$) in the
above inequality we deduce that relation \eq{11} holds true. The
proof of Lemma \ref{l1} is complete.  \endproof

\begin{lemma}\label{l2}
There exists $u\in E\setminus\{0\}$ such that $\frac{J(u)}{I(u)}=\lambda_1$.
\end{lemma}
\proof Let $\{u_n\}\subset E\setminus\{0\}$ be a minimizing
sequence for $\lambda_1$, that is,
\begin{equation}\label{stea}
\lim_{n\rightarrow\infty}\frac{J(u_n)}{I(u_n)}=\lambda_1>0.
\end{equation}
By relation \eq{10} it is clear that $\{u_n\}$ is bounded in $E$. Since $E$ is reflexive it follows that
there exists $u\in E$ such that $u_n$ converges weakly to $u$ in $E$. On the other hand, similar arguments
as those used in the proof of Lemma 3.4 in \cite{RoyalSoc} show that the functional $J$ is weakly lower
semi-continuous. Thus, we find
\begin{equation}\label{14}
\liminf_{n\rightarrow\infty}J(u_n)\geq J(u).
\end{equation}
By relation \eq{3} it follows that
$E$ is compactly embedded in $L^{q(x)}(\Omega)$. Thus, $u_n$ converges strongly in $L^{q(x)}(\Omega)$. Then, by
relation \eq{L6} it follows that
\begin{equation}\label{15}
\lim_{n\rightarrow\infty}I(u_n)=I(u).
\end{equation}
Relations \eq{14} and \eq{15} imply that if $u\not\equiv 0$ then
$$\frac{J(u)}{I(u)}=\lambda_1.$$
Thus, in order to conclude that the lemma holds true it is enough
to show that $u$ is not trivial. Assume by contradiction the
contrary. Then $u_n$ converges weakly to $0$ in $E$ and strongly
in $L^{q(x)}(\Omega)$. In other words, we will have
\begin{equation}\label{victor}
\lim_{n\rightarrow\infty}I(u_n)=0.
\end{equation}
Letting $\epsilon\in(0,\lambda_1)$ be fixed by relation \eq{stea} we deduce that for $n$ large enough we have
$$|J(u_n)-\lambda_1I(u_n)|<\epsilon I(u_n),$$
or
$$(\lambda_1-\epsilon)I(u_n)<J(u_n)<(\lambda_1+\epsilon)I(u_n).$$
Passing to the limit in the above inequalities and taking into account that relation \eq{victor} holds true we find
$$\lim_{n\rightarrow\infty}J(u_n)=0.$$
That fact combined with relation \eq{L6} implies that actually $u_n$ converges strongly to $0$ in $E$, i.e.
$\lim_{n\rightarrow\infty}\|u_n\|=0$. By this information and relation \eq{11} we get
 $$\lim_{n\rightarrow\infty}\frac{J(u_n)}{I(u_n)}=\infty,$$
 and this is a contradiction. Thus, $u\not\equiv 0$.
 The proof of Lemma \ref{l2} is complete.  \endproof

 By Lemma \ref{l2} we conclude that there exists $u\in E\setminus\{0\}$ such that
\begin{equation}\label{16}
\frac{J(u)}{I(u)}=\lambda_1=\inf_{w\in E\setminus\{0\}}\frac{J(w)}{I(w)}.
\end{equation}
Then, for any $v\in E$ we have
$$\frac{d}{d\epsilon}\frac{J(u+\epsilon v)}{I(u+\epsilon v)}\left|_{\epsilon=0}=0\right..$$
A simple computation yields
\begin{equation}\label{17}
\int_{\Omega}(|\nabla u|^{p_1(x)-2}+|\nabla u|^{p_2(x)-2})\nabla u\nabla v\;dx\cdot I(u)-
J(u)\cdot\int_\Omega|u|^{q(x)-2}uv\;dx=0,\;\;\;\forall\;v\in E.
\end{equation}
Relation \eq{17} combined with the fact that $J(u)=\lambda_1 I(u)$
and $I(u)\neq 0$ implies the fact that $\lambda_1$ is an eigenvalue
of problem \eq{1}. Thus, step 2 is verified.
\medskip

\noindent$\bullet$ {\sc Step 3.} We show that any
$\lambda\in(\lambda_1,\infty)$ is an eigenvalue of problem \eq{1}.

Let $\lambda\in(\lambda_1,\infty)$ be arbitrary but fixed. Define $T_\lambda:E\rightarrow\RR$ by
$$T_\lambda(u)=J(u)-\lambda I(u).$$
Clearly, $T_\lambda\in C^1(E,\RR)$ with
$$\langle T_\lambda^{'}(u),v\rangle=\langle J^{'}(u),v\rangle-\lambda\langle I^{'}(u),v\rangle,\;\;\;\forall\;u\in E.$$
Thus, $\lambda$ is an eigenvalue of problem \eq{1} if and only if there exists $u_\lambda\in E\setminus\{0\}$ a
critical point of $T_\lambda$.

With similar arguments as in the proof of relation \eq{10} we can
show that $T_\lambda$ is coercive, i.e.
$\lim_{\|u\|\rightarrow\infty}T_\lambda(u)=\infty$. On the other
hand, as we have already remarked, similar arguments as those used
in the proof of  Lemma 3.4 in \cite{RoyalSoc} show that the
functional $T_\lambda$ is weakly lower semi-continuous. These two
facts enable us to apply Theorem 1.2 in \cite{S} in order to prove
that there exists $u_\lambda\in E$ a global minimum point of
$T_\lambda$ and thus, a critical point of $T_\lambda$. In order to
conclude that step 4 holds true it is enough to show that
$u_\lambda$ is not trivial. Indeed, since $\lambda_1=\inf_{u\in
E\setminus\{0\}}\frac{J(u)}{I(u)}$ and $\lambda>\lambda_1$ it
follows that there exists $v_\lambda\in E$ such that
$$J(v_\lambda)<\lambda I(v_\lambda),$$
or
$$T_\lambda(v_\lambda)<0.$$
Thus,
$$\inf_{E}T_\lambda<0$$
and we conclude that $u_\lambda$ is a nontrivial critical point of
$T_\lambda$, or $\lambda$ is an eigenvalue of problem \eq{1}. Thus,
step 3 is verified.
\medskip

\noindent$\bullet$ {\sc Step 4.} Any $\lambda\in(0,\lambda_0)$,
where $\lambda_0$ is given by \eq{lambda0}, is not an eigenvalue
of problem \eq{1}.

Indeed, assuming by contradiction that there exists $\lambda\in(0,\lambda_0)$ an eigenvalue of problem \eq{1}
it follows that there exists $u_\lambda\in E\setminus\{0\}$ such that
$$\langle J^{'}(u_\lambda),v\rangle=\lambda\langle I^{'}(u_\lambda),v\rangle,\;\;\;\forall\;v\in E.$$
Thus, for $v=u_\lambda$ we find
$$\langle J^{'}(u_\lambda),u_\lambda\rangle=\lambda\langle I^{'}(u_\lambda),u_\lambda\rangle,$$
that is,
$$J_1(u_\lambda)=\lambda I_1(u_\lambda).$$
The fact that $u_\lambda\in E\setminus\{0\}$ assures that
$I_1(u_\lambda)>0$. Since $\lambda<\lambda_0$, the above
information yields
$$J_1(u_\lambda)\geq\lambda_0 I_1(u_\lambda)>\lambda I_1(u_\lambda)=J_1(u_\lambda).$$
Clearly, the above inequalities lead to a contradiction.  Thus, step
4 is verified.
\medskip

By steps 2, 3 and 4 we deduce that $\lambda_0\leq\lambda_1$. The
proof of Theorem \ref{t1} is now complete. \endproof

\smallskip
\noindent{\bf Remark 1.} At this stage we are not able to deduce
whether $\lambda_0=\lambda_1$ or $\lambda_0<\lambda_1$. In the
latter case an interesting question concerns the existence of
eigenvalues of problem \eq{1} in the interval
$[\lambda_0,\lambda_1)$. We propose to the reader the study of
these open problems.


\begin{thebibliography}{99}
{\footnotesize

\bibitem{acerbi2} E. Acerbi and G. Mingione, Gradient estimates for the $p(x)$-Laplacean system,
{\it J.~Reine Angew. Math.} {\bf 584} (2005), 117-148.

\bibitem{chen} Y. Chen, S. Levine and R. Rao, Functionals with $p(x)$-growth in
image processing, Duquesne
University, Department of Mathematics and Computer Science Technical Report 2004-01, available
at {\tt www.mathcs.duq.edu/\~{}sel/CLR05SIAPfinal.pdf}.

\bibitem{D} L. Diening, {\it Theoretical and Numerical Results for
Electrorheological Fluids}, Ph.D. thesis, University of Frieburg,
Germany, 2002.

\bibitem{edm} D. E. Edmunds, J. Lang, and A. Nekvinda, On $L^{p(x)}$
norms, {\it Proc. Roy. Soc. London Ser.~A} {\bf 455} (1999), 219-225.

\bibitem{edm2} D. E. Edmunds and J. R\'akosn\'{\i}k, Density of smooth
functions in $W^{k,p(x)}(\Omega)$, {\it Proc. Roy. Soc. London Ser.~A}
{\bf 437} (1992), 229-236.

\bibitem{edm3} D. E. Edmunds and J. R\'akosn\'{\i}k, Sobolev embedding
with variable exponent, {\it Studia Math.} {\bf 143} (2000),
267-293.

\bibitem{FZh} X. L. Fan and Q. H. Zhang, Existence of solutions for
$p(x)$-Laplacian Dirichlet problem, {\it Nonlinear Anal} {\bf 52}
(2003), 1843-1852.

\bibitem{FZZ} X. Fan, Q. Zhang and D. Zhao, Eigenvalues of $p(x)$-Laplacian
Dirichlet problem, {\it J. Math. Anal. Appl.} {\bf 302} (2005), 306-317.


\bibitem{hal} T. C. Halsey, Electrorheological fluids, {\it Science}
{\bf 258} (1992), 761-766.


\bibitem{KR} O. Kov\'a\v cik and J. R\'akosn\'{\i}k, On spaces
$L^{p(x)}$ and
$W^{1,p(x)}$, {\it Czechoslovak Math. J.} {\bf 41} (1991), 592-618.

\bibitem{RoyalSoc} M. Mih\u ailescu and V. R\u adulescu,
A multiplicity result for a nonlinear degenerate problem arising in
the theory of electrorheological fluids, {\it Proceedings of the Royal Society A: Mathematical, Physical and Engineering Sciences}
{\bf 462} (2006), 2625-2641.

\bibitem{mihradproc} M. Mih\u ailescu and V. R\u adulescu, On a nonhomogeneous quasilinear
eigenvalue problem in Sobolev spaces with variable exponent, {\it
Proc. Amer. Math. Soc.} {\bf 135} (2007), 2929-2937.

\bibitem{mihradjmaa} M. Mih\u ailescu and V. R\u adulescu, Existence and multiplicity of solutions for quasilinear
nonhomogeneous problems: an Orlicz-Sobolev space setting, {\it
J.~Math. Anal. Appl.} {\bf 330} (2007), 416-432.

\bibitem{M} J. Musielak, {\it Orlicz Spaces and Modular  Spaces},
Lecture Notes in Mathematics, Vol. 1034, Springer, Berlin, 1983.

\bibitem{R} M. Ruzicka, {\it Electrorheological Fluids: Modeling
and Mathematical Theory}, Springer-Verlag, Berlin, 2002.

\bibitem{samko} S. Samko and B. Vakulov, Weighted Sobolev theorem with variable exponent
for spatial and spherical potential operators, {\it J. Math. Anal. Appl.} {\bf 310} (2005), 229-246.

\bibitem{S} M. Struwe, {\it Variational Methods: Applications to
Nonlinear Partial Differential Equations and Hamiltonian Systems},
Springer, Heidelberg, 1996.

\bibitem{Z1} V. Zhikov, Averaging of functionals in the calculus of
variations and elasticity, {\it Math. USSR Izv.} {\bf 29} (1987),
33-66.

}
\end{thebibliography}
\end{document}